\documentclass[a4paper,11pt]{article}

\usepackage{geometry,amsmath,amssymb,latexsym,enumerate,xypic}
\def\id{\operatorname{id}}

\def\biglabels{\def\labelstyle{\textstyle}}

\newtheorem{example}{Example}[section]
\newtheorem{Def}[example]{Definition}
\newtheorem{prop}[example]{Proposition}
\newtheorem{Theo}[example]{Theorem}

\newtheorem{cor}[example]{Corollary}

\def\geq{\geqslant}

\def\Ext{\operatorname{Ext}}
\def\Aut{\operatorname{Aut}}

\newenvironment{prf}{{\noindent \bf Proof} }{\hfill $\Box$
\def\biglabels{\def{\labeltyle{\textstyle}}}

\mbox{}}

\begin{document}

\title{\large\bf COVERING GROUPS OF NON-CONNECTED TOPOLOGICAL
GROUPS REVISITED}          \author{Ronald Brown and Osman Mucuk
\\ School of Mathematics \\ University of Wales \\
           Bangor, Gwynedd LL57 1UT, U.K.}

\date{Received 11 March 1993 ,  Revision 31 May 1993}

\maketitle {\em Published in Math. Proc. Camb. Phil. Soc},  115
(1994) 97-110.

This version: 18 Dec. 2006, with some new references.
\section*{Introduction}

     All spaces are assumed to be locally path connected and
semi-locally 1-connected.   Let  $X$   be  a  connected
topological group with identity   $e$, and let
$p:\tilde{X}\rightarrow X$ be the universal cover of the
underlying space of $X$.    It follows  easily  from classical
properties  of lifting maps to covering spaces that for any
point  $\tilde{e}$  in  $\tilde{X}$  with  $p\tilde{e} = e$ there
is a unique structure  of  topological  group on $\tilde{X}$ such
that  $\tilde{e}$  is the identity and $p\colon
\tilde{X}\rightarrow X$ is a morphism of groups. We say that the
structure of topological group on  $X$   {\em lifts} to
$\tilde{X}$.

It is less generally appreciated that this result fails for  the
non-connected  case. The set $\pi_{0}X$ of path components  of
$X$  forms a non-trivial group which acts on the abelian group
$\pi_{1}(X,e)$ via conjugation in $X$. R.L. Taylor~\cite{Ta2}
showed that the topological group $X$ determines  an {\em
obstruction class}  $k_X$  in $H^3(\pi_{0}X,\pi_{1}(X,e))$, and
that the  vanishing  of  $k_X$ is a  necessary  and  sufficient
condition  for  the  lifting of  the  topological group structure
on  $X$  to a  universal covering  so  that  the  projection is a
morphism.  Further, recent work, for  example
Huebschmann~\cite{Hu2}, shows there are geometric applications of
the non-connected case.

The purpose of this paper is to prove generalisations of this
result on coverings of topological groups using modern work  on
coverings of groupoids (see for example, Higgins~\cite{Hi},
Brown~\cite{Br2}),  {\em via} the following scheme. We first use
the fact that covering spaces of space $X$ are equivalent to
covering morphisms of the fundamental groupoid $\pi _1X$ (section
1). This extends easily to the group case: if $X$ is a
topological group, then the fundamental groupoid inherits a group
structure making it what is called a {\em group-groupoid}, i.e.
a  group  object  in  the  category of groupoids; then
topological group coverings of $X$ are equivalent to
group-groupoid coverings of $\pi_1X$ (Proposition 2.3).

The next input is the equivalence between group-groupoids,  and
crossed modules (Brown and Spencer~\cite{Br-Sp}).   Here a
crossed module is a morphism $\mu :M\rightarrow P$ of groups
together with an action of the group $P$ on the group $M$, with
two axioms satisfied. It is easy to translate notions from
covering morphisms of group-groupoids  to corresponding notions
for crossed modules (Proposition 4.2).

The existence of simply connected covering groups of a
topological group now translates  to the existence of extensions
of groups of ``the type of a given crossed module'' (Definition
5.1).   This generalisation of the classical extension theory is
due to Taylor \cite{Ta1} and Dedecker ~\cite{De}. We formulate a
corresponding notion of abstract kernels (Theorem 5.2), analogous
to that due to Eilenberg-Mac Lane ~\cite{Mac}. This leads to our
main result, Theorem 5.4,  which determines when a morphism
$\theta : \Phi\rightarrow \pi_0X$ of groups is realised by a
covering morphism $p : \tilde{X}\rightarrow X$ of topological
groups such that $\tilde{X}  $  is simply connected with $\pi_0X
$ isomorphic to $\Phi$.  We deduce that any topological  group $X$
admits a simply connected covering group covering all the
components of $X$ (Corollary 5.6).  According to comments in
\cite{Ta2}, results of this type were known to Taylor.

Our proof of Theorem 5.2 uses methods of crossed complexes, as in
Brown and Higgins \cite{Br-Hi1}. This seems the natural setting
for these  results, since crossed complexes contain information
on resolutions and on crossed modules.  The exposition is
analogous to that  of Berrick ~\cite{Be} for  the ordinary
theory  of  extensions,  in that fibrations are used, but  in
the  algebraic context of crossed complexes.  A direct account of
a special case of these results, in the context of Lie groupoids,
is given by Mackenzie in~\cite{Ma}, and this account could also
be adapted to the general case.

Section 6 deals with coverings other than simply connected ones.

     The results of this paper formed part of Part I of
Mucuk~\cite{Mu}.

\section{Groupoids and coverings}

The  main  tool  is  the  equivalence  between  covering  maps
of   a  topological space  $X$   and covering morphisms of the
fundamental groupoid   $\pi_{1}X$  of  $X$.  Our main reference
for groupoids and this result  is Brown~\cite{Br2}  but  we adopt
the following notations and conventions.

A topological space $X$ is called {\em simply connected\/} if
each loop in $X$ is contractible in $X$, and $X$ is called 1-{\em
connected\/}  if it is connected and simply connected. A map
$f:X\rightarrow Y$ is called $\pi _0${\em -proper} if $\pi _0
(f)$ is a bijection.

If $X$  is a topological  space,  the  category  $ TCov/X $ of
covering spaces of  $X$  is the full subcategory of the slice
category   $Top/X$  of spaces over    $X$  in which the objects
are  the  covering  maps.   It  is  standard  that  if $h\colon
Y\rightarrow Z$  is a map in   $TCov/X$, i.e. is a  map over
$X$,   then  $h$   is  a  covering map.  Further,  if $f\colon
Y\rightarrow X$  is a covering map such that  $Y$  is simply
connected,  then  for any other cover  $g\colon Z\rightarrow X$,
there is a covering map  $h\colon Y\rightarrow Z$   over   $X$.
This is summarised by saying that  $Y$  covers any other cover
of   $X$, and a covering map with this property is called {\em
universal}.  A necessary and sufficient condition for this is
that $Y$ be simply-connected.

For a groupoid  $G$, we write  $O_G$  for the set of objects of
$G$,        and   $G$  for the set of arrows, or elements.  We
write $s, t\colon G\rightarrow O_G$   for   the  source  and
target maps.   The product $g \circ h$  is defined if and only
if  $tg = sh$.    The  identity at  $x\in O_G$  is written
$1_x$.  The inverse of an  element    $g$   is written
$g^{-1}$.   \par

The category of groupoids and morphisms of groupoids is written
$Gd$.\par

For  $x\in O_G$  we denote the {\em star}    $\{g\in G\mid sg =
x\}$   of  $x$   by  $G^x$, and the {\em costar}   $\{g\in G\mid
tg = x\}$  of  $x$  by  $G_x$, and write $G^x_y$ for $G^x\cap
G_y$.  The object  group  at  $x$  is  $G(x) = G^x_x$. An
element of some $G^x_x$ is called a {\em loop\/} of $G$. \par

We say  $G$  is {\em transitive} (resp. 1-{\em transitive},
{\em simply  transitive})  if  for all  $x, y\in O_G$, $G(x,y)$
is non-empty (resp. is a singleton,  has not more  than one
element).    The transitive  component  of  an object $x$   of
$G$   is  the largest transitive subgroupoid of $G$ with $x$ as
an object, and is written  $C(G,x)$.  The set of  transitive
components  of   $G$  is   written  $\pi_{0}G$.   A morphism $p$
of groupoids is called $\pi _0${\em -proper} if $\pi _0 (p)$ is
a bijection.  \par

Covering morphisms and universal covering groupoids of a
groupoid   are defined in Brown~\cite{Br1} (see  also
Higgins~\cite{Hi},   Brown~\cite{Br2})  as  follows: \par

 Let  $p\colon H\rightarrow G$  be a morphism of groupoids.
Then   $p$   is  called  a {\em covering morphism} if for each
$x\in O_H$, the restriction  $G^{x} \rightarrow G^{px}$     of
$p$   is bijective.  The covering morphism $p$ is called {\em
regular} if for all objects $x$ of $G$ and all $g\in G(x)$ the
elements of $p^{-1}(g)$ are all or none of them loops. This is
equivalent to the condition that for all objects $y$ of $H$, the
subgroup $pH(y) $ of $G(py)$ is a normal subgroup \cite{Br2}.
\par

If $G$  is a groupoid, the category $GdCov/G$ of coverings of
$G$ is the full subcategory of the slice category $Gd/G$ of
groupoids over $G$ in which the objects are the covering
morphisms.\par

A covering morphism  $p\colon H\rightarrow G$  is called {\em
universal} if  $H$  covers every covering of   $G$,  i.e.  if
for  every  covering  morphism   $a\colon A\rightarrow G$ there
is   a   morphism    of   groupoids   $a'\colon H\rightarrow A$
such that   $aa' = p$   (and  hence   $a'$ is  also  a  covering
morphism).   It is common to consider universal covering
morphisms which are $\pi _0$-proper.  \par

 We recall the following standard  result   (Brown~\cite{Br2},
Chapter  9), which summarises the theory of covering  spaces.

\begin{prop}  For any space  $X$,     the fundamental groupoid
functor   defines  an equivalence of categories

\[ \pi_1: TCov/X  \rightarrow GdCov/(\pi_{1}X). \]

  \end{prop}

One crucial step in the  proof  of  this  equivalence  is  the
result  (Brown~\cite{Br2}, 9.5.5) that if   $q\colon
H\rightarrow \pi_{1}X$   is   a  covering  morphism  of
groupoids, then there is a topology on  $O_H$  such that
$O_{q}\colon O_{H}\rightarrow X$   is  a  covering map,  and
there  is  an  isomorphism
$\alpha\colon\pi_{1}O_{H}\rightarrow H$   such  that   $q\alpha
= \pi_{1}(O_q)$.  This result,  which  translates  the  usual
covering  space  theory into a more base-point free context,
yields the inverse equivalence.  \par

We also remark that the  universal cover of $X$ at  $x\in X$  is
given by the target map   $(\pi_{1}X)^{x}\rightarrow X$with the
subspace topology from a topology on $\pi _1 X$ .   \par

 Recall that an {\em action\/} of a groupoid  $G$  on sets {\em
via}  $ w$   consists of a function $w\colon A\rightarrow O_G$,
where  $A$  is a set, and  an  assignment  to  each   $g\in
G(x,y)$   of  a  function    $g_{\sharp}\colon
w^{-1}(x)\rightarrow  w^{-1}(y)$,   written   $a\mapsto  a\circ
g$, satisfying the usual rules  for  an  action,   namely
$a\circ 1 = a$,  $a\circ (g\circ h) = (a\circ g)\circ h$ when
defined.  A {\em morphism}   $f\colon (A,w)\rightarrow (A',w')$
of such actions is a  function  $f\colon A\rightarrow A'$  such
that  $w'f = w$  and    $f(a\circ g) = (fa)\circ g$   whenever
$a\circ g$  is defined.  This gives a category  $Act(G)$  of
actions of   $G$  on  sets.   For such an action, the {\em
action groupoid} $A\rtimes G$  is defined to   have  object  set
$A$, arrows the pairs  $(a,g)$  such that $ w(a) = sg$, source
and   target  maps  $s(a,g) = a$, $t(a,g) = a\circ g$, and
composition                  \[           (a,g)\circ (b,h) =
(a,g\circ h) \]                          whenever  $b = a\circ
g$.  The projection    $q\colon A\rtimes G\rightarrow  G$,
$(a,g)\mapsto g$,  is  a  covering morphism of groupoids, and
the functor sending an action  to  this  covering morphism gives
an equivalence of  categories    $Act(G) \rightarrow GdCov/G$.
(See for example Brown~\cite{Br2}.) \par

Let  $x$  be an object of the transitive groupoid  $G$, and let
$N(x)$  be  a subgroup of the object group $G(x)$.  Then  $G$
acts on  the  set   $A$    of  cosets  $N(x)\circ g$  for  $g\in
G^x$, {\em via}  the map   $N(x)\circ g\mapsto tg$.  So we can
form  the corresponding covering morphism  $p\colon H\rightarrow
G$,  where    $H = A\rtimes G$,  and  the object  $\tilde{x} =
N(x)$  of  $H$  satisfies   $p(H(\tilde{x})) = N(x)$.  This
construction  yields  an  equivalence  of  categories  between
the  lattice  ${\cal L}G(x)$    of  subgroups of  $G(x)$  and
the category of pointed transitive coverings of  $G,x$.   \par

Suppose   further   that    $a\in G^x_y$,        $ N(y) =
a^{-1}\circ N(x)\circ a$,     and  $q\colon K\rightarrow G$  is
the covering of $G$  determined as above  by    $N(y)$,   with
$\tilde{y}\in O_K$  satisfying  $q[K(\tilde{y})] = N(y)$.  Then
there is   a  unique  isomorphism   $h\colon H\rightarrow K$
such that  $qh = p$  and $h\tilde{z} =\tilde{y}$.  That is,
conjugate  subgroups  of a transitive groupoid $G$ determine
isomorphic coverings, and we  obtain  an  equivalence of
categories between  the  lattice  of  conjugacy  classes  of
subgroups of  $G$  and the isomorphism classes  of  transitive
coverings  of   $G$.   \par

If $G$  is not transitive then $\pi _0$-proper coverings may be
constructed by  working  on each transitive component.  We
choose a transversal for the set   $I =\pi_{0}G$   of components
of  $G$, i.e. an object  $\tau_i$  for each component  $G_i$
of   $G$,  choose a subgroup  $N(\tau_i)\subseteq G(\tau_i)$,
and get a covering   $\tilde{G_i}\rightarrow G_{i}$   for  each
component  $G_i$  of $G$.  The disjoint union of these coverings
is a  covering   $p\colon \tilde{G}\rightarrow G$, which is
universal if and only if  all the  $N(\tau_i)$  are  trivial
groups.


\section{Group-groupoids and covering morphisms}

The notion of group-groupoid, and the first parts of
Propositions 2.1 and 2.3 below are taken from Brown-Spencer
\cite{Br-Sp}, although the term used there is $\cal G$-groupoid.
\par

By a {\em group-groupoid\/} we mean a groupoid  $G$  with a
morphism  of  groupoids  $G\times G\rightarrow G$, $(g,h)\mapsto
gh$,   yielding a group  structure internal to the  category  of
groupoids.   Since  the  multiplication  is  a  morphism   of
groupoids, we obtain the {\em interchange  law},  that
$(a\circ g)(b\circ h) = (ab)\circ (gh)$,  for all $g,h,a,b\in G$
such that   $a\circ g$   and   $b\circ h$   are  defined.   If
the  identity for the group structure on  $O_G$  is written $e$,
then  $1_e$   is  the  identity for the group structure on the
arrows.  The group  inverse  of  an  arrow  $g$ is  written
$\bar{g}$.   Then   $g\mapsto \bar{g}$  is  a   morphism
$G\rightarrow G$   of  groupoids.   \par

It is a standard consequence of the interchange law  that  the
groupoid composition  in a group-groupoid can be recovered from
the group law, as shown in the first  part  of  the following
proposition.

\begin{prop}  Let  G  be a group-groupoid, and suppose  $a\circ
b$  is  defined  in  G , where $a\in G(x,y)$.  Then  $a\circ b =
a{\bar{1}_y}b$.   If further   $g\in G(e)$, then    \[
a\circ (1_{y}g)\circ a^{-1}  = 1_{x}g, \]     and     \[
ag\bar{a} = 1_{x}g{\bar{1}_x}.  \]                    Further,
G(e)  is abelian.

 \end{prop}

\begin{prf} Suppose  $ta = y$.  Then
\begin{align*}
a\circ b & = ((a{\bar{1}_y})1_{y})\circ (1_{e}b) \\
& = ((a{\bar{1}_y})\circ 1_{e})(1_{y}\circ b) \\           & =
a{\bar{1}_y}b .
\end{align*}
Further \begin{alignat*}{2}
 a\circ (1_{y}g)\circ
a^{-1} &  =  a\circ ((1_{y}g)\circ (a^{-1}1_{e}))& & = a\circ
((1_{y}\circ a^{-1})(g\circ 1_{e})) \\   &= a\circ (a^{-1}g)& &=
(a1_{e})\circ (a^{-1}g) \\                                &=
(a\circ a^{-1})(1_{e}\circ g) &&= 1_{x}g . \\ \intertext{On the
other hand}
 a\circ (1_{y}g)\circ a^{-1}
&  =  (a \bar{1}_y 1_y g)\circ a ^{-1} && =  (ag)\circ (1_y
\bar{a} 1_x) \\ & =  ag\bar{1}_y 1_y  \bar{a} 1_x && =  ag\bar{a}
1_x .
\end{alignat*}
Hence  $ag\bar{a}   = 1_x g\bar{1}_x$.        That  $a\circ g =
g = g\circ a$  for $a,g\in G(e)$  is immediate.
\end{prf}
 \begin{cor} Let  $N(e)$  be a subgroup of   $G(e)$,  and  let
$N$ be  the  family of subsets  $N(x) = 1_{x}N(e)$  for  all
$x\in O_G$.    Then   $N(x)$   is  a  normal subgroupoid of  $G$.
In particular, all the object groups of  $G$   are  isomorphic,
and are abelian.
\end{cor}
\begin{prf} That  $N(x)$  is a subgroup follows from  \[
 1_{x}(b\circ a) = (1_{x}\circ 1_{x})(b\circ a) = (1_{x}b)\circ
(1_{x}a) , \]                  for $b,a\in N(x)$.  The normality
follows from the second formula   of   the  Proposition, on
taking  $g\in N(e)$.  It is  immediate  that  all  the  object
groups are isomorphic.
\end{prf}
This result  implies  that  all  coverings  of  a group-groupoid
are  regular.   It  also  shows  that  a  choice $\tau$   of
transversal  for   the components of a group-groupoid   $G$
induces  an  equivalence  between   the category ${\cal L}G(e)$
of subgroups of  $G(e)$  under inclusion and the category  of
isomorphism classes of $\pi_{0}$-proper coverings of  $G$.   \par

We now consider coverings in the category of group-groupoids.
\par

A {\em  morphism}  of  group-groupoids  is  a  morphism  of the
underlying  groupoids which preserves the group structure.
Then  group-groupoids  and  morphisms of them form a category
which we will denote by $GpGd$. Let $G$ be a group-groupoid.
Then $GpGdCov/G$ denotes the full subcategory of the slice
category $GpGd/G$ whose objects are group-groupoids
$p:H\rightarrow G$ over $G$ such that $p$ is a covering morphism
of the underlying groupoids.

     We can now translate Proposition 1.1 to this situation.

\begin{prop}  Let  X  be a  topological  group. Then  the
fundamental  groupoid   $\pi_{1}X$  is a group-groupoid with
group structure induced by that  of   X .  Further,  the
fundamental  groupoid  functor   $\pi_{1}$  gives  an
equivalence  from  the category $GpTCov/X$ to the category
$GpGd/\pi_{1}X$. \end{prop}

\begin{prf} We show that the inverse equivalence of Proposition
1.1  determines  an inverse equivalence in this case also. \par
Suppose then that  $q\colon H\rightarrow \pi_{1}X$  is a
morphism  of group-groupoids such  that the underlying groupoid
morphism is a covering morphism.   Then  there  is a topology on
$\tilde{X} = O_H$   and  an  isomorphism
$\alpha\colon\pi_{1}\tilde{X}\rightarrow H$   such  that   $p =
O_{q}\colon\tilde{X}\rightarrow X$  is a covering map and
$q\alpha = \pi_1(p)$.  The  group  structure  on $H$  transports
via  $\alpha$  to a morphism of groupoids \[
\tilde{m}\colon\pi_{1}\tilde{X}\times\pi_{1}\tilde{X}\rightarrow
\pi_{1}\tilde{X} \]              \noindent                such
that   $\pi_{1}(p)\circ \tilde{m} = m\circ
(\pi_{1}(p)\times\pi_{1}(p))$,    where    $m$    is   the
group  multiplication on $X$, and clearly  $\tilde{m}$  is a
group structure on    $\pi_{1}\tilde{X}$.   By  9.5.5 of
Brown~\cite{Br2}, $\tilde{m}$  induces a continuous map on
$\tilde{X}$.   This  gives  the  multiplication on  $\tilde{X}$.
The fact that this is  a  group structure  follows  from the
fact  that $\tilde{m}$  is a group structure.  \end{prf}


 \section{Actions of group-groupoids on groups}     In this
section we relate group-groupoid  covering  morphisms to a
notion of action of a group-groupoid on a  group.   The  results
are a special case of results of Section 1 of Brown and
Mackenzie~\cite{Br-Ma}, and are included here for completeness.
\par

 Let  $G$  be a group-groupoid.    An {\em action } of the
group-groupoid $G$  on a group  $A$       {\em via} $w$
consists  of a morphism  $w\colon A\rightarrow G$  from the
group  $A$  to the  underlying group of  $O_G$   and an action
of the groupoid  $G$  on the  underlying  set   $A$    via   $w$
such that the following interchange law holds: \[       (a\circ
g)(b\circ h) = (ab)\circ (gh)  \]
\noindent whenever both sides are defined.  A {\em morphism}
$f\colon (A,w)\rightarrow (A',w')$  of such  operations is a
morphism  $f\colon A\rightarrow A'$  of  groups  and  of   the
underlying operations of  $G$.  This gives a category
$GpGdAct(G)$.  For an  action  of  $G$  on the group  $A$  via
$w$, the action  groupoid  $A\rtimes G$   is defined.  It
inherits a group structure by \[        (a,g)(c,k) = (ac,gk) .
\]                       It is easily checked that   $A\rtimes
G$   is  then  a  group-groupoid,  and  the  projection
$p\colon A\rtimes G\rightarrow G$  is an object of the category
$GpGd/G$.  By  means of this construction, one obtains the
following, which is  a  special  case of Theorem 1.7 of Brown
and Mackenzie~\cite{Br-Ma} which considers  the  case  of
actions of Lie double groupoids.

\begin{prop}  The  categories GpGdCov/G and  GpGdAct(G) are
equivalent.
             \end{prop}


\section{Group-groupoids and crossed modules}     A {\em crossed
module}  $(M,P,\mu)$  is defined in      Whitehead~\cite{Wh} to
consist  of two groups $M$  and $P$  together with a
homomorphism   $\mu\colon M\rightarrow P$, and an  action of
$P$  on  $M$  on the right,  written  $(m,p)\mapsto m^p$,
such  that  the  following conditions are satisfied:  \par  CM1)
$\mu(mp) = p^{-1}(\mu m)p$ \par CM2)  $n^{\mu m}   = m^{-1}nm$
\par \noindent for all $m,n\in M$  and  $p\in P$. \par
Standard examples of crossed modules are:

\begin{enumerate}[(i)]
  \item the inclusion $M\rightarrow P$   of a normal subgroup,
   \item the zero morphism  $M\rightarrow P$   when    $M$
is  a   P-module,
  \item the inner  automorphism  map
$\chi _M: M\rightarrow \Aut r
M$   for   any  group $M$,
\item a morphism  $M\rightarrow P$  of groups which  is
surjective  and   has  central  kernel,
\item the free crossed $P$-module  $C(w)\rightarrow
P$  arising from  a function   $w\colon R\rightarrow P$  (see
Brown and Huebschmann~\cite{Br-Hu}),
\item the  induced morphism
$\pi_{1}(F,x)\rightarrow \pi_{1}(E,x)$  of fundamental groups
for   any  fibration  of  spaces   $F\rightarrow E\rightarrow B$.
\end{enumerate}

Standard consequences of the  axioms   (see for example
\cite{Br-Hu}) are  that   $\mu M$   is  a  normal subgroup of
$P$,  that  $Ker\mu$  is central in $M$,  and that   $\mu M$
acts  trivially on  $Ker\mu$ which thereby becomes a module over
$Coker\mu$. \par

A {\em morphism}  $(f,g)\colon (M,P,\mu)\rightarrow (N,Q,\nu)$
of      crossed modules consists of  group morphisms  $f\colon
M\rightarrow N$  and  $g\colon P\rightarrow Q$   such that
$g\mu =\nu f$   and   $f$   is an operator homomorphism, that
is,   $f(m^p) = {fm}^{(gp)}$      for    $m\in M$   and   $p\in
P$.  So  crossed  modules  and  morphisms  of  them,  with the
obvious  composition of morphisms  $(f',g')(f,g) = (f'f,g'g)$,
form a category,  which  we write   $CrsM$.  \par

The following theorem was found by Verdier in 1965, but not
published,  and found independently by Brown and
Spencer~\cite{Br-Sp}.  We give  a   sketch  of  the proof, since
we need some of its detail.

\begin{Theo}   The category  $GpGd$  of group-groupoids is
equivalent  to  the  category   $CrsM$  of crossed modules.  If
a group-groupoid  G  has associated  crossed module  $(M,P,\mu)$
then the underlying groupoid of  G  is   transitive  (resp.
simply  transitive,  1-transitive)  if  and  only  if    $\mu$
is   an  epimorphism (resp. a monomorphism, isomorphism).
Further, the  group    $\pi_{0}G$   is  $Coker\mu$. \end{Theo}

{\bf Sketch Proof:} A functor $\delta\colon GpGd\rightarrow
CrsM$  is  defined as  follows.   For  a  group-groupoid  $G$
we let  $\delta (G)$  be the crossed module  $(M,P,\mu)$   where
$P$   is  the  group    $O_G$   of  objects of  $G$; $M$   is
the costar  $G_e$  of  $G$  at the identity  $e$ of the group
$O_G$;  $\mu\colon M\rightarrow P$  is the restriction of the
source map   $s$ ;  the  group structures on  $M$ and  $P$  are
induced by that on  $G$ ; and   $P$    acts  on   $M$   by
$m^p$ =${\bar{1}_p}m1_p$   for   $p\in P$   and  $m\in M$.
The  results   on  transitivity follow immediately.\par
Conversely define a functor  $\beta\colon CrsM\rightarrow GpGd$
in the following  way.  For a crossed module   $(M,P,\mu)$,
$\beta (M,P,\mu)$   is  the   group-groupoid  whose  object set
(group) is $P$  and whose group of  arrows  is  the  semi-direct
product   $P\ltimes M$  with the standard group structure  \[
(p,m)(q,n)=(pq, m^{q} n). \]                         The source
and target maps   $s, t$   are  defined  to  be    $s(p,m) = p$
and   $t(p,m) = p(\mu m)$,  while the composition of arrows is
given by            \[ (p,m)\circ (q,n) = (p,mn)   \]
        whenever  $p(\mu m) = q$.  \hfill   $\Box$

If  $X$  is a topological group with identity  $e$, then the
fundamental  groupoid  $\pi_{1}X$  becomes a group-groupoid, the
associated crossed   module is   $t\colon (\pi_{1}X)^{e}
\rightarrow X$  (Brown and Spencer~\cite{Br-Sp}), and
$(\pi_{1}X)^{e}$  has a topology making it  the universal  cover
based  at  $e$   of   the  path component of  $e$.   \par

 It is easy to obtain results for morphisms of group-groupoids
corresponding  to  Theorem 4.1, as follows.

 \begin{prop} Let  $f\colon H\rightarrow G$  be a morphism of
group-groupoids and  let  $(f_{1},f_{2})\colon
(N,Q,\nu)\rightarrow (M,P,\mu)$   be   the   morphism   of
crossed   modules  corresponding to  $f$  as in Theorem 4.1.
Then, on underlying groupoids,    $f$   is a covering morphism if
and only if    $f_{1}\colon  N\rightarrow M$   is  an
isomorphism.   Further,  $f$  is a universal covering morphism if
and only  if   $f_1$    is  an  isomorphism,   $\nu$    is   a
monomorphism,   and   the   induced    morphism
$Coker\nu\rightarrow Coker\mu$  is an isomorphism. \end{prop}

We therefore define a morphism  $(f_{1},f_{2})$  of crossed
modules      as in the  proposition to be a {\em covering
morphism} if  $f_1$  is  an  isomorphism,   and  so  obtain a
category   $CrsMCov/(M\rightarrow P)$  of coverings of
$M\rightarrow P$   as  a  full  category of the slice category
$CrsM/(M\rightarrow P)$.

\begin{cor} The category $GTCov/X$  of topological group
coverings  of  a  topological group  $X$  is  equivalent  to the
category   $CrsMCov/((\pi_{1}X)^{e}\rightarrow X)$   of  crossed
module  coverings  of   $(\pi_{1}X)^{e}\rightarrow{X}$.

 \end{cor}


\section{Extensions, crossed modules and \newline  cohomology}

We now recall the notion of an extension of groups of  the  type of
a  crossed module, due to Taylor~\cite{Ta1} and Dedecker~\cite{De}.
See also \cite{B-P}.

\begin{Def} {\em Let   $\cal{M}$   denote  the  crossed  module
$\mu \colon M\rightarrow P$.    An  {\it extension $(i,p,\sigma)$
of type  $\cal{M}$ }  of the group   $M$  by the group   $\Phi$ is
first  an exact sequence of groups  $$ \xymatrix{ 1  \ar [r] & M
\ar [r]^{i} & E\ar [r]^{p}&\Phi\ar [r]&1 } $$

\noindent so that  $E$  operates on  $M$  by conjugation, and
$i\colon M\rightarrow  E$   is  hence  a  crossed module.
Second, there is given  a morphism of crossed modules

$$ \xymatrix {  1  \ar  [r]  & M \ar @{=} [d]\ar  [r] ^{i} & E \ar  [d] ^{\sigma}\\        & M
\ar  [r] _{\mu} &P   }   $$ i.e.  $\sigma i =\mu$  and  $m^{e} =
m^{\sigma e}$, for all $m\in M$,  $e\in  E$.} \end{Def}

Two such extensions of type  $\cal{M}$ $$\biglabels \xymatrix {  1
\ar [r] & M  \ar  [r] ^{i} & E\ar  [r] ^{p}&\Phi\ar  [r] &1,\\ 1 \ar
[r] & M  \ar  [r] ^{i'} & E'\ar  [r] ^{p'}&\Phi\ar  [r] &1,
               }   $$ are said to be {\em equivalent} if  there is
a  morphism  of  exact  sequences
$$ \biglabels \xymatrix {  1
\ar [r] & M  \ar @{=} [d]\ar  [r] ^{i} & E\ar  [r] ^{p} \ar
[d]^\phi&\Phi\ar @{=}
[d]\ar  [r] &1,\\
1 \ar [r] & M  \ar  [r] ^{i'} & E'\ar  [r] ^{p'}&\Phi\ar  [r] &1,
               }  $$ \noindent such that
$\sigma'\phi = \sigma$. Of course in this case  $\phi$ is an
isomorphism, by the 5-lemma, and hence equivalence of extensions is
an equivalence relation. Denote by  $\Ext_{\cal M} (\Phi, M)$ the
set of equivalence classes of all extensions of  type $\cal{M}$ of
$M$ by $\Phi$.

 An extension of   $M$   by $\Phi$   of type   $\cal{M}$
determines a  morphism   $\theta \colon \Phi \rightarrow Q$, where
$Q = Coker\mu$,  which  is dependent only   on   the equivalence
class of the extension, and  $\theta $ is here called the {\em
abstract ${\cal M}$-kernel}  of the extension.  The set of extension
classes with  a  given  abstract ${\cal M}$-kernel $\theta$  is
written $\Ext_{(\cal{M},\theta)}(G,M)$.   \par

 The usual  theory  of  extensions  of   a group $M$ by a group
$\Phi$ considers extensions of  the  type  of the crossed module
$\chi_{M} \colon M\rightarrow \Aut M$.   The advantages of replacing
this by a general crossed module are first  that the group   $\Aut
M$ is  not  a functor of $M$, so that the relevant cohomology
theory in  terms  of    $\chi_{M}$ appears to have no coefficient
morphisms, and second, that the more general  case occurs
geometrically, as in \cite{Ta2} and in this paper.   \par We now
show there is an obstruction to realizability, analogous to the
classical result of Eilenberg-Mac Lane (\cite[Ch.V, Prop.8.3]{Mac}).
The cohomology groups $H^{\ast}_{\theta}(\Phi,A)$ referred to here
are defined later.

\begin{Theo} Let ${\cal M}$ be the crossed module $\mu : M
\rightarrow P$ with $A = Ker\mu,\,Q = Coker \mu$. Let $\theta :
\Phi \rightarrow Q$ be an abstract ${\cal M}$-kernel . Then
there is an {\em obstruction class} $k({\cal M},\theta ) \in
H^3_{\theta } (\Phi ,A)$ whose vanishing is necessary and
sufficient for there to exist an extension of M by $\Phi $ of
type ${\cal M}$ with abstract $\cal{M}$-kernel $\theta $.
Further, if the obstruction class is zero, then the equivalence
classes of such extensions are bijective with $H^2_{\theta }
(\Phi ,A)$.  \end{Theo}

We give an exposition of a proof of this theorem using the methods
of crossed complexes as given for example in Brown and
Higgins~\cite{Br-Hi1} or~\cite{Br-Hi3}. The point is that crossed
complexes allow for methods analogous to those of chain complexes as
in standard homological algebra, but including non-abelian
information of the type given by crossed modules. The obstruction
result arises from an exact sequence of a fibration of crossed
complexes. This allows us to give a proof analogous to that given
for the classical case using topological methods by Berrick
in~\cite{Be}. A direct proof may also be given by extending the
methods of Mackenzie \cite{Ma} to more general crossed modules than
$M\rightarrow \Aut(M) $

        We assume the definition of crossed complex  as
given  for example     in  Brown and Higgins~\cite{Br-Hi1}
or~\cite{Br-Hi3}, and in particular  the notion  of  pointed
morphism.  Recall that a {\em reduced} crossed  complex has a
single vertex.  A {\em homotopy} $h\colon f\simeq g$  of pointed
morphisms   $f,g\colon C\rightarrow D$   of crossed complexes is
a family of functions   $h_{i} \colon C_{i}\rightarrow D_{i+1}$
such  that   \par  i)  $h_{1}\colon C_{1}\rightarrow D_{2}$  is
a derivation over  $g_1$,   that is, \[    h_{1}(x + y) =
{h_{1}(x)}^{gy}   + h_{1}(y), \]       where  $g(y) = g_{1}(y)$,
for  $x,y\in C_1$.          \par ii) For $n\geq 2$,
$h_{n}\colon C_{n}\rightarrow D_{n+1}$  is an  operator morphism
over  $g_1$, that is,    \[          h_{n}(x^{a}  + y) =
(h_{n}x)^{ga}   + h_{n}(y), \]      where  $ga = g_{1}a$.
   \par iii)  If  $x\in C_{1}$, then        \[        gx =
fx\delta h_{1}x  . \]                                \par iv) If
$n\geq 2$  and  $c\in C_{n}$, then        \[     gx =
fxh_{n-1}\delta x - \delta h_{n}x . \]
\par            We will also use the morphism crossed complex
$CRS_{\ast}(C,D)$        defined  in  Brown and
Higgins~\cite{Br-Hi2} whose elements in dimension  0  are   the
pointed  morphisms  $C\rightarrow D$, in dimension  1  are  the
homotopies,   and  in  higher  dimensions are the ``higher
homotopies".  \par     A  crossed module  $\mu \colon M\rightarrow
P$  can also be extended by trivial groups to give a crossed complex
\[ \cdots \rightarrow 1 \rightarrow \cdots 1 \rightarrow
1\rightarrow M\rightarrow P. \]  Denote this crossed complex again
by  $\cal M$.

Let  $\Phi$  be  a  group.   We  write  $C\Phi$   for  the
standard  crossed  resolution of  $\Phi$.  This is defined in
Huebschmann~\cite{Hu} and   shown  in Brown and
Higgins~\cite{Br-Hi1} to be the fundamental crossed complex of
the (Kan) simplicial set,  $Nerv(\Phi)$, the nerve of the group
$\Phi$.   \par

Write  $[C\Phi, \cal{M}]$  for the set of pointed homotopy
classes  of       morphisms  $C\Phi\rightarrow M$.

\begin{Theo}  There is a bijection                            \[
[C\Phi, {\cal M}]\cong \Ext_{{\cal M}} (\Phi,M).\]
     \end{Theo}

The proof is given in Brown and Higgins~\cite{Br-Hi1}.  The key
point      is  that $C_{1}\Phi$  is the free group on elements
$[g]$,      $g\in \Phi$, $C_{2}\Phi$  is the free crossed
$C_{1}\Phi$-module on   $\delta \colon \Phi \times \Phi
\rightarrow C_{1}\Phi$,  where   $\delta (g,h) =
[g][h][gh]^{-1}$, and,  for  $i\geq 3$, $C_{i}\Phi$  is the
free  $\Phi$-module  on   $[g_{1},\ldots ,g_{i}]$, for
$g_{1},\ldots, g_{i}\in \Phi$. Further because  of the form of
the  boundary  morphism  $\delta \colon C_{3}\Phi\rightarrow
C_{2}\Phi$, a morphism   $C\Phi\rightarrow \cal{M}$  is
equivalent to  a  factor  set (with values in ${\cal M}$), and a
homotopy  of  morphisms  is  essentially  an  equivalence of
factor sets.   \par

Recall that   $\cal{M}$   is  the  crossed  module        $\mu
\colon M\rightarrow P$,   and    $A = Ker\mu$,  $Q = Coker\mu$. Let
$\xi \cal{M}$,  $\zeta \cal{M}$    denote  the  crossed complexes in
the following diagram of morphisms of crossed complexes

$$\biglabels\xymatrix @C=2.5pc{  \ldots\ar  [r]&1 \ar  [d]   \ar [r] &M   \ar  [d] _{\id} \ar
[r]& P     \ar  [d] _{\id}&& \ar  [d]  {\cal M}\\
\ldots \ar [r] & A \ar  [d] _{\id} \ar [r]   & M  \ar  [d]\ar [r] &
P \ar  [d] && \ar  [d]_{q} \xi
\cal{M}\\
\ldots \ar [r] & A \ar [r]   & 1   \ar [r] &Q &&\zeta \cal{M}
 }$$

\noindent where  $q$  is determined by the quotient morphism $P
\rightarrow Q$.  Since  $q$  is an  epimorphism   in  each
dimension, it is also  a  fibration  of  crossed  complexes  and
therefore,  since  $C\Phi$  is free, the induced morphism of
morphism complexes \[  q_{\ast}\colon CRS_{\ast}(C\Phi,\xi {\cal
M})\rightarrow CRS_{\ast}(C\Phi,\zeta {\cal M}) \]
is also a  fibration  of  crossed  complexes   (Brown  and
Higgins~\cite{Br-Hi3}, Prop.6.2).   Since  $C\Phi$  is free and
$\xi M$  is acyclic, there is an identification   \[
\pi_{0}CRS_{\ast}(C\Phi,\xi {\cal M})\cong Hom(\Phi,Q) .\]
           Further, each morphism $\theta \colon \Phi\rightarrow
Q$  determines an  action of  $\Phi$  on  $A$   and  so a cohomology
group  $H_{\theta}^{3}(\Phi,A)$.  Then
$\pi_{0}CRS_{\ast}(C\Phi,\zeta \cal{M})$  is the union of  all these
cohomology groups for all such  $\theta$.  The function
$\pi_{0}(q_{\ast})$   takes  a  morphism  $\theta $  to a cohomology
class         \[     k({\cal M},\theta)\in H_{\theta}^{3} (\Phi,A) ,
\]                           called the {\em obstruction class} of
$({\cal M},\theta)$.    If $k\colon C\Phi\rightarrow \xi \cal{M}$
is a realisation of $\theta$, then  $qk$  represents
$k(\cal{M},\theta )$.   If this  class  is  0,  then there is a
homotopy  $h\colon qk \simeq l$,  say,  where    $l_{1} = qk_{1}$,
$l_3 = 0$.   Hence $k_{3} = h_{2}\delta$.   So  there  is  a
homotopy    $k\simeq k'$   where   ${k'}_{1} = k_1$,  ${k'}_{2} =
k_{2} - \delta h_{2}$, ${k'}_{3} = 0$.   \par     Let  $F$  be the
fibre of $q_{\ast}$  over  $l$.  Then       $\pi_{0}F$  may  be
identified  with the set $[C\Phi ,\cal{M}]$ of homotopy classes of
morphisms $C\Phi \rightarrow {\cal M}$, and so with the classes  of
extensions  of  $A$  by  $\Phi$  of type $\cal{M}$.  The exact
sequence of  the fibration    $q_{\ast}$ with  fibre  $F$  yields,
given the above identifications, the exact sequence \[ 0\rightarrow
H_{\theta}^{2}(\Phi ,A)\rightarrow  \Ext_{\cal M}(\Phi
,M)\rightarrow Hom(\Phi ,Q)\rightarrow  H_{\theta}^{3}(\Phi ,A)
\;\;\;  (\star)
           \] \noindent where the three right hand terms have
base points the class of the split extension, the morphism $\theta$,
and zero respectively. The obstruction part of Theorem 5.2 follows
immediately.  The standard theory of the exact sequence of a
fibration of crossed complexes~\cite{Ho} also yields that the group
$H_{\theta}^{2}(\Phi ,A)$  operates on $\Ext_{\cal{M}}(\Phi ,M)$
so  that  the  classes  of extensions of type  $\cal{M}$  with
abstract kernel  $\theta$ are given  by this group.  This completes
the proof of Theorem 5.2. $\Box $ \par

We can translate Theorem 5.2 to the following.

\begin{Theo} Let  $X$  be a topological group. Let $\Phi$ be a
group, and   let $\theta\colon\Phi\rightarrow \pi_{0}X$  be a
morphism  of groups.  Then there is  a  covering morphism
$p\colon \tilde{X}\rightarrow X$  of topological  groups   and
an  isomorphism   $\alpha\colon\pi_{0}\tilde{X}\rightarrow\Phi$
such that   $\theta\alpha = \pi_{0}(p)$  and  $\tilde{X}$  is
simply  connected if and only if the obstruction class   \[
k({\cal M},\theta)\in H_{\theta}^{3} (\Phi ,\pi_{1}(X,e)) \]
               is  zero,  where $\cal{M}$   is  the  associated
crossed  module    $(\pi_{1}X)^{e}\rightarrow X$.  Further, the
isomorphism classes of such coverings are bijective with
$H^2_{\theta}(\Phi ,A)$.  \end{Theo}

\begin{prf} We write  $\mu\colon M\rightarrow X$  for
$\cal{M}$.  If the  obstruction  class  is  zero  then there is
an extension   $1\rightarrow M\rightarrow E\rightarrow
\Phi\rightarrow 1$  of  type    $\cal{M}$,   and  the  crossed
module  $M\rightarrow E$  corresponds to a simply transitive
group-groupoid   $\tilde{G}$.  The morphism from   $M\rightarrow
E$  to  $\cal{M}$ yields a  covering morphism of
group-groupoids   $\tilde G\rightarrow \pi_{1}X$.    Hence  we
obtain  the  required  covering  space  $\tilde{X} =
Ob(\tilde{G})$.   The  converse   follows  from  Theorem 5.2, as
does the classification of these coverings.   \end{prf}

If $\cal M$ is an arbitrary crossed module with cokernel $Q$,
and one takes $\Phi = Q$ and $\theta = {\rm id}$ in 5.2,  then
the class      $k({\cal M}, {\rm id})\in H^{3}(Q,A)$, where
the  action  of $Q$  on  $A$  is the given one, is called the
{\em obstruction  class}   $k(\cal{M})$   of  the crossed module
$\cal{M}$.  As a consequence of Theorem 5.4  we  recover   the
result of Taylor~\cite{Ta2}.

\begin{cor}  Let  $X$  be a (possibly disconnected) topological
group and  let    $p\colon\tilde{X}\rightarrow X$   be  a  $\pi
_0$-proper  universal covering.  Then the group structure of
$X$  lifts to  $\tilde{X}$   such that   $\tilde{X}$  is a
topological group and  $p$  is a morphism of   topological
groups  if  and only if the obstruction class   $k({\cal M})\in
H^{3}(\pi_{0}X, \pi_{1}(X,e))$   is zero. \end{cor}

 We remark that this obstruction class is shown in      Brown
and  Spencer~\cite{Br-Sp}, to  be  the  first
$k$-invariant  of  the  classifying  space  of  the  topological
group  $X$.  \par

The following result is referred to in \cite{Ta2}.

\begin{cor} Let $X$ be a (possibly disconnected) topological
group. Then there exists a simply connected covering group
$p:\tilde{X}\rightarrow X$ of $X$  such that $\pi _0 p$ is
surjective.

\end{cor}

\begin{prf} It is enough to choose an epimorphism $\theta : \Phi
\rightarrow \pi _0 X$ such that the induced morphism on
cohomology

\[ \theta ^{\ast} :H^{3}(\pi_{0}X, \pi_{1}(X,e))\rightarrow
H_{\theta}^{3} (\Phi ,\pi_{1}(X,e)) \] is trivial. This can be
done with $\Phi $ a free group.

\end{prf}

Of course, there is no uniqueness result for this simply
connected cover. \par

In the next section, we generalise Theorem 5.4 to  a  wider
class  of  coverings.   \section{General coverings of
topological groups}

We now deal with other coverings than simply connected ones, as does
Taylor in \cite{Ta2} for the proper case. \par We first recall two
basic constructions which will be used later.  The first essentially
gives  the  usual  forward  coefficient morphism   in cohomology.
\begin{prop} Let $\mu\colon M\rightarrow P$  be a crossed module
with  $A = Ker\mu$  and $Q = Coker\mu$.  Let $\phi\colon
A\rightarrow B$  be a morphism of Q-modules.  Then there is  a
crossed module  ${\mu}'\colon M'\rightarrow P$  and a morphism of
exact  sequences $$\biglabels \xymatrix
{  0\ar  [r] &A\ar  [r] ^{i}\ar  [d] _{\phi}&M\ar  [r] ^{\mu}\ar  [d] _{{\phi}'}&P\ar  [r]\ar  [d]  &Q\ar  [d] \ar  [r] &1 \\
0\ar  [r] &B\ar  [r] ^{j}&M'\ar  [r] ^{{\mu}'}&P\ar  [r] &Q\ar [r]
&1
 }   $$
such that $(\phi ', {\rm id})$ is a morphism of crossed modules.

 \end{prop}

\begin{prf}  The proof is easy on taking  $M' = (B \times M)/C$,
where    $C = (\phi ,i)(A)$,  and defining  ${\mu}'$  by
$[b,m]\mapsto \mu m$,  ${\phi}'$   by    $m\mapsto [m,1]$,
where   $[b,m]$   denotes the class of  $(b,m)$  in  $M'$.
\end{prf}

 \begin{prop} Let   $\cal{M}$   be  the  crossed  module
$\mu\colon M\rightarrow P$,   let   $Q = Coker \mu$, and let
$\theta\colon\Phi\rightarrow Q$  be an abstract kernel.  Then
\[         k({\cal M},\theta) = k({\cal N},{\rm id}) \]
            where  $\cal N$   is the crossed module $\nu\colon
M\rightarrow P\times_{Q}\Phi$,   $m\mapsto (m,1)$.   Further, there
is a bijection \[  \Ext_{({\cal M},\theta)}(\Phi ,M)\cong
\Ext_{({\cal N},{\rm id})}(\Phi,M) .\]
\end{prop}

 \begin{prf} This follows from the morphism of exact sequences

$$\biglabels\xymatrix {
0\ar  [r] &A\ar  [d] _{\id}\ar  [r]
^{i}&M\ar  [d] _{\id}\ar  [r] ^-{\nu}&P\times_{Q}\Phi^{\rule{0em}{1.2ex}}\ar  [r] \ar  [d]&\Phi\ar  [d] _{\theta} \ar  [r] &1 \\
0\ar  [r] &A\ar  [r] ^{i}&M\ar  [r] ^{\mu}&P\ar  [r] &Q\ar [r] &1
 }   $$ \end{prf}

Now we can give the following theorem.

 \begin{Theo} Let X be a topological group, let
$\theta\colon\Phi\rightarrow\pi_{0}X$ be a morphism of groups,
and let  N be a $\pi_{0}X$-invariant subgroup of
$\pi_{1}(X,e)$. Then  there is a covering morphism $p\colon
\tilde{X}\rightarrow X$ of  topological groups and an
isomorphism  $\alpha\colon\pi_{0}\tilde{X}\rightarrow \Phi$ such
that  $\theta\alpha=\pi_{0}(p)$ and
$p(\pi_{1}(\tilde{X},\tilde{e}))=N$  if and only if the
obstruction class \[ k({\cal M},\theta)\in H_{\theta}^{3}(\Phi,
\pi_{1}(X,e)), \] where $\cal{M}$  is the associated crossed
module   $(\pi_{1}X)^{e}\rightarrow X$, is mapped to zero by the
morphism induced by the coefficient morphism   \[
\pi_{1}(X,e)\rightarrow  (\pi_{1}(X,e))/N .   \]
\end{Theo}

 \begin{prf} Write  the  crossed   module    $\cal{M}$    as
$\mu\colon M\rightarrow P$,    and   let   $Q = Coker\mu$, $A =
Ker\mu$.  Suppose that there is such a  covering   morphism  of
topological groups and isomorphism  $\alpha$  as given in  the
theorem.    Let the crossed module  $\cal{N}$  associated to
$\tilde{X}$  be written as $\nu : \tilde{M} \rightarrow E$, so that
$Ker \nu = N$. Then $\cal N$ maps to $\cal M$   as part of the
following diagram  $$\biglabels \xymatrix {  0\ar  [r] &N\ar  [d]
_{i} \ar  [r] &\tilde{M}\ar  [d] _{\cong}
\ar  [r] ^{\nu}&E     \ar  [d] _{\sigma}        \ar  [r] &\Phi \ar  [d] _{\theta}\ar  [r] &1  \\
 0\ar [r] &A\ar  [r]   &M \ar  [r] ^{\mu}
&P \ar [r]  &Q    \ar  [r] &1 }
$$ where $N$ is now a $\Phi$-module via  $\theta$.    Let  $\cal
{M}'$  and ${\cal M} \rightarrow {\cal M}'$ be  the  crossed
module and morphism of crossed modules constructed from the
quotient mapping $A\rightarrow A/N$  as   in  Proposition  6.1.
Let  $k\colon C\Phi\rightarrow \xi {\cal N}$  be a realisation of
the identity morphism on   $\Phi$.   Then  the  composite
$C\Phi\rightarrow \xi{\cal N}\rightarrow \zeta {\cal N}$
realises   $k({\cal N},{\rm id})$.   Clearly  the  composition \[
C\Phi\rightarrow \xi {\cal N}\rightarrow \zeta {\cal N}\rightarrow
\zeta {\cal M}\rightarrow \zeta {\cal M'} \]
realises the zero class in  $ H_{\theta }^{3}(\Phi,A/N)$, as
required.
\par     Suppose conversely that  $k({\cal M},\theta)$  maps to
zero  in         $ H_{\theta }^{3}(\Phi,A/N)$.   Again,  let $\cal
M'$  be the crossed module constructed in Proposition 6.1, with
morphism $\phi ' :M \rightarrow M'$.  Then, by assumption, the
obstruction  class  $k( {\cal M'},\theta )$  is zero, and so there
is an extension of type $\cal{M}'$ and with abstract kernel $\theta$
$$ \xymatrix {  & &1\ar  [r] &M'\ar  [r] ^{i'}&E\ar [r] &\Phi\ar
[r] &1&    }   $$
   \noindent    It is easy to check that      $\nu =
i'{\phi}'\colon M\rightarrow E$  becomes a crossed  module  when
$E$ acts on $M$ via $\sigma$, and that  $Ker\nu = Ker{\phi}' = N$.
Hence  we  have  the  following  morphism  of  exact sequences
$$ \xymatrix {
0\ar  [r] &N\ar  [d] \ar  [r] &M\ar  [d]\ar  [r] ^{\nu}&E\ar  [d] _{\sigma}\ar  [r] &\Phi\ar  [d] _{\theta}\ar  [r] &1& \\
0\ar  [r] &A\ar  [r] &M\ar  [r] ^{\mu}&P\ar  [r] &Q\ar  [r] &1 }
$$
The morphism of crossed modules this includes can be realised by a
covering morphism of group-groupoids  and  so of  topological groups
as required. \end{prf}

\begin{example} {\em We mention some nice examples of Taylor
\cite{Ta2}. He shows there are exactly three non-isomorphic
topological group extensions of $SO(2)$ by $\bf{Z}_2$, namely
the direct sum of the two groups, the orthogonal group $O(2)$,
and finally the multiplicative group of all quaternions
$a+bi+cj+dk$, of norm 1, such that $(a^2 + b^2)(c^2 + b^2)=0$.
Other examples of non-connected coverings of topological groups
are given in section 8 of \cite{Ta2}. } \end{example} \par

This completes our account of the theory of covering groups of
topological groups. \par

     Of  course  these   theorems  on   spaces   have
analogues   for  group-groupoids which we leave the reader to state.
\section*{Acknowledgements}     We would like to thank: Kirill
Mackenzie and Johannes Huebschmann  for  helpful comments and
conversations, and a referee for many helpful comments.  The second
author would like to  thank  the Turkish Government for support
during his studies at Bangor.

 \end{document}